\documentclass{svmult}
\usepackage[bottom]{footmisc}
\usepackage{bbm,enumerate,url}
\newcommand{\kla}{\left ( }
\newcommand{\mer}{\right ) }
\newcommand{\noo}{\left \|}
\newcommand{\rrm}{\right \|}
\newcommand{\bet}{\left |}
\newcommand{\rag}{\right |}
\newcommand{\equa}{\begin{eqnarray*}}
\newcommand{\tion}{\end{eqnarray*}}
\newcommand{\bes}{\widetilde{\mathbbm{B}}}
\newcommand{\Om}{\Omega}
\newcommand{\ftn}{{\cal F}}
\newcommand{\sptext}[3]{\hspace{#1 em}\mbox{#2}\hspace{#3 em}}
\renewcommand{\P}{\mathbbm{P}}
\newcommand{\Q}{\mathbbm{Q}}
\newcommand{\R}{\mathbbm{R}}
\renewcommand{\E}{\mathbbm{E}}
\renewcommand{\D}{\mathbbm{D}}
\newcommand{\B}{\mathbbm{B}}
\def\SDE{{\sc (SDE)}}
\def\SDEs{{\sc (SDE) }}
\def\GBM{{\sc (GBM)}}
\def\GBMs{{\sc (GBM) }}

\begin{document}

\title{Fractional smoothness and applications in Finance}
\author{Stefan Geiss \and  Emmanuel Gobet}
\institute{Department of Mathematics \\
           University of Innsbruck \\
           Technikerstra\ss e 13/7 \\
           A-6020 Innsbruck \\
           Austria \\
           \texttt{stefan.geiss@uibk.ac.at}
\and       Laboratoire Jean Kuntzmann \\
           Universit\'e de Grenoble and CNRS \\
           BP 53, 38041 Grenoble Cedex 9 \\
           France \\
           \texttt{emmanuel.gobet@imag.fr}}

\maketitle

\begin{abstract}
This overview article concerns the notion of fractional smoothness of 
random variables of the form $g(X_T)$, where $X=(X_t)_{t\in [0,T]}$ is 
a certain diffusion process. We review the connection to the real 
interpolation theory, give examples and applications of this 
concept. The applications in stochastic finance mainly concern
the analysis of discrete time hedging errors. We close the review 
by indicating some further developments.
\end{abstract}

\noindent


\section{Introduction}
\label{sec:introduction}

From the practitioners one learns that hedging an option which payoff is discontinuous is more difficult compared to the case of smooth payoffs: this feature appears for instance for digital options or barrier options (we refer the reader to \cite{tale:92} among others). On the one hand, for such options the number of assets (i.e. the delta) to incorporate in the hedging portfolio is unbounded, and it may become larger and larger when one gets close to the singularity (i.e. the maturity and the strike for digital options, or the trigger level for barrier options). On the other hand, the numerical estimation of this delta becomes less and less accurate, leading to global stability issues. These heuristic observations are the starting point for
 deeper mathematical investigations about the concept of irregular payoffs, in order to formalize it and to quantify the
 payoff irregularity (with the notion of {\em fractional smoothness}). In the current contribution, we aim to give an
 overview of this concept and some applications in stochastic finance. Actually, the applications go beyond the
 financial framework and more generally, they concern the theory of stochastic differential equations and their
 approximations.

\paragraph{The discrete time hedging error as an important application.}
Since most of the results presented here are applied to the aforementioned example of hedging possibly irregular options, 
we start with a brief presentation of this problem, in order to emphasize the issues to handle and to raise some natural 
questions. Take for instance an European-style option exercised at maturity $T>0$, with a payoff of the form $h(S_T)$ 
where $S_t:=[S^1_t,\cdots,S^d_t]$ denotes the price of a $d$-dimensional underlying asset at time
$0\le t \le T$. Sometimes, we will use the notation $X^i_t=\log(S^i_t)$ for the log-asset, and $g(x)=h(e^{x_1},\cdots,e^{x_d})$ 
for the payoff in the logarithmic variables. In what follows, we assume a Markovian dynamics without jumps for the 
asset (solution to a SDE defined below), we suppose that the interest rate is equal to 0 (to simplify
the presentation) and that the market is complete (for details about this standard framework, see \cite{kara:shre:98}). 
Thus, under some regularity assumptions, the payoff $h(S_T)$ can be replicated perfectly by a continuous time strategy, 
where $\delta^S_t=\nabla_x H(t,S_t)$ defines the vector of number of assets to hold at time $t$. Here, $H$ is the fair price 
of the option, that is
\[ H(t,x)=\E_\Q(h(S_T)|S_t=x) \]
where $\Q$ is the (unique) risk-neutral measure. In practice, only discrete-time hedging is possible at some times  $\tau=(t_i)_{i=0}^n$ with $0=t_0 < t_1 < \cdots < t_{n-1} < t_n = T$. Thus, at time $t\in[0,T]$ the option seller is left with the tracking error
\begin{eqnarray}
      C_t(h(S_T),\tau)
& = & H(t,S_t)-H(0,S_0)-\sum_{i} \nabla_x H(t_i,S_{t_i})\cdot (S_{t_{i+1}\land t}-S_{t_{i}\land t})
      \nonumber \\
& = & \int_0^t (\nabla_x H(s,S_s)-\nabla_x H(\phi(s),S_{\phi(s)}))\cdot dS_s\label{eq:track:error}
\end{eqnarray}
with $\phi(s)=t_i$ when $t_i< s\le t_{i+1}$.
We expect the tracking error (\ref{eq:track:error}) to converge to 0 as the number $n$ of 
re-balancing dates goes to infinity.
With the above formulation (\ref{eq:track:error}), the tracking error is naturally associated to 
the problem of approximation of a stochastic integral using piece-wise constant integrated processes. But the delta process  $(\nabla_x H(s,S_s))_{s\in [0,T)}$ may exhibit very different behaviors from payoff to payoff: if the payoff is smooth enough, then the delta might be {\em bounded} as time goes
to maturity,  while an irregular payoff usually yields an {\em exploding} delta as $s\rightarrow T$. This gives rise to the first question.
\begin{itemize}
\item [{\bf (Q1)}] \quad 
      Is there an intrinsic way to relate the growth rate (as $s\rightarrow T$) of the         derivatives of $H$ to the irregularity of the payoff $h$?
      \smallskip

      The answer will be {\em yes} via the notion of {\em fractional smoothness} introduced below, 
      see Theorems \ref{thm:gm_1} and \ref{thm:gm_2}.
\end{itemize}
The estimation of stochastic integrals is usually performed with $L_2$-norms, but in our financial setting, both measures $\P$ and $\Q$ can be considered. For practitioners, errors under the historical probability $\P$ are presumably more relevant, while the mathematical treatment under the risk-neutral measure $\Q$ is simpler in our context
(because the tracking error process (\ref{eq:track:error}) is a $\Q$-local martingale).
\begin{itemize}
\item [{\bf (Q2)}] \quad
      Is the definition of fractional smoothness affected by the choice of a specific
      measure? 
      Do the $L_2$-convergence rates depend on the choice of the probability measures $\P$ or $\Q$?
      \smallskip 

      In the context we consider the answer concerning the fractional smoothness 
      is usually {\em no} in the sense of the comments after 
      Theorem \ref{thm:gm_2}. Concerning the approximation rates the same is checked for examples so far
      (see the remarks after Theorem \ref{thm:gobet-temam-binary}).
\end{itemize}
Beyond the approach to measure tracking errors in $L_2$, we could alternatively identify the weak limit of the re-normalized tracking error.
\begin{itemize}
\item [{\bf (Q3)}] \quad
      Do the weak convergence rates coincide with those in the $L_2$ sense?
      \smallskip

      The answer is {\em not necessarily}, as there are counter-examples in which the convergence in $L_2$ 
      and in distribution hold at different rates, see Section \ref{section:appli}.
\end{itemize}
Finally, through an efficient choice of re-balancing dates $\tau$, one can expect to reduce tracking errors and improve the risk management of options.
\begin{itemize}
\item [{\bf (Q4)}] \quad
      Which time-nets $\tau=(t_i)_{i=0}^n$ lead to optimal convergence rates? 
      And how to relate them to the  {\em fractional smoothness} of the payoff?
      \smallskip

      As answer we get that according to the index of fractional smoothness of the payoff, one
      can define explicitly re-balancing times achieving the optimal convergence rates, see
      Section \ref{section:appli}.
\end{itemize}
These preliminary questions serve as references for the reader when reading the next sections. 

\paragraph{Organization of the paper.} First, we define
the probabilistic framework and the assumptions used throughout this work. Then in Section
\ref{section:2}, we define the {\em fractional smoothness} and provide basic properties: we
choose a presentation that is quite illuminating regarding the previous preliminary questions.
In Section  \ref{sec:interpolation}, we take another view on {\em fractional smoothness} using
the interpolation theory. In Section \ref{sec:examples}, we consider examples of terminal conditions 
and identify their fractional smoothness. Then, in Section \ref{section:appli}, we go back to the
 analysis of discrete time hedging errors and state the main results. We close by further developments 
 and applications of the {\em fractional smoothness} 
 in Section \ref{sct:devpt}.

\paragraph{Assumptions.}  Let us define the probabilistic setting used in the following. We fix a $d$-dimensional 
Brownian motion $W=(W_t)_{t\in [0,T]}$ 
defined on a complete probability space $(\Omega,\ftn_T,\P)$ and we let 
$(\ftn_t)_{t\in [0,T]}$ be the augmentation of the natural filtration of $W$. 
The log-asset $X$ is the solution of the $d$-dimensional forward diffusion
\[ X_t = x_0 + \int_0^t b(s,X_s)ds + \int_0^t \sigma(s,X_s) dW_s. \]
To state the results, we mainly consider two types of assumptions:

\begin{itemize}
\item [{\bf (SDE)}] \hspace*{2.5em}
      $d\ge 1$ and $b,\sigma \in C_b^\infty([0,T]\times \R^d)$ and $\sigma \sigma^* \ge \delta I_{\R^d}$ 
      for some $\delta>0$. 
\item [{\bf (GBM)}] \hspace*{2.4em}
      $d= 1$ and $X_t=\ln(S_t) = W_t - (t/2)$.
\end{itemize}
The smoothness conditions in \SDEs are too strong, they are chosen to simplify the presentation. Whenever useful to simplify even more, we may consider the very simple case of the geometric Brownian motion \GBMs (here, the asset is a martingale, meaning that $\P=\Q$). The reader is referred to the corresponding original papers for the possible weaker conditions. 
\medskip

In the following $|\cdot|$ stands for the Euclidean norm and $A\sim_c B$ for $A/c \le B \le cA$ if
$c\ge 1$ and $A,B\ge 0$. Expectations and conditional expectations under $\P$ are simply denoted by $\E(.)$ and $\E(.|{\cal F}_t)$, while under $\Q$, we indicate explicitly the dependency w.r.t. the probability measure by writing $\E_\Q(.)$ and $\E_\Q(.|{\cal F}_t)$.


\section{Definition of fractional smoothness and basic properties}
\label{section:2}
Fractional smoothness on the Wiener space can be defined in various ways,
see \cite{wata:93,hirs:99}.
Our approach is motivated by the questions discussed in Section \ref{sec:introduction}.
Since we consider only random variables of the form $Z=g(X_T)=h(S_T)$ (a function of the process
at maturity $T$), the time $T$ plays a specific role in our definition. It would be necessary to
modify our definition for more general dependencies like 
$Z=g(X_{t_1},\cdots, X_{t_n})$, see \cite{geis:geis:gobe:09}.
\medskip

\begin{definition}\label{def:Besov_X} \rm
Assume that $Z\in L_2(\P)$.
\begin{enumerate}[{\rm (i)}]
\item For $0<\theta \le 1$ we let $Z \in \bes_{2,\infty}^\theta$ provided that, for
      all $0\le t <T$,
      \[     \| Z - \E(Z|\ftn_t) \|_{L_2(\P)}
               \le c (T-t)^\frac{\theta}{2}.\]
\item For $0<\theta<1$ we let $Z\in \bes_{2,2}^\theta$ provided that
      \[   \int_0^T (T-t)^{-1-\theta} \| Z - \E(Z|\ftn_t) \|_{L_2(\P)}
^2 dt 
         < \infty.\]
\end{enumerate}
\end{definition}
\smallskip
The spaces $\bes_{2,q}^\theta$ above will always be obtained by the conditional expectation and the
$L_2$-norm under the measure $\P$. Therefore we omit the dependency on $\P$ in the notation.

The following properties follow straight from the definition:
\smallskip

\begin{proposition}
For $0<\theta < \eta < 1$ and $p,q\in\{ 2,\infty \}$ we have that
      \[            \bes_{2,\infty}^1 
         \subseteq  \bes_{2,p}^\eta
         \subseteq  \bes_{2,q}^\theta
         \sptext{1}{and}{1}
                    \bes_{2,2}^\theta
         \subseteq  \bes_{2,\infty}^\theta.\]
\end{proposition}
\bigskip
Given a bounded
\footnote{Here again, the boundedness assumptions on $g$ can be weakened and we refer to 
          the original papers.}  
measurable $g:\R^d\to \R$ and $(t,x)\in [0,T)\times \R^d$, we let
\equa
u(t,x)    &:= & \E (g(X_T) | X_t=x), \\
D^2u(t,x) &:= & \kla \frac{\partial^2 u}{\partial x_i \partial x_j}(t,x) 
                \mer_{i,j=1}^d.
\tion
The following equivalences are useful to exploit properties of
$\bes_{2,2}^\theta$ and $\bes_{2,\infty}^\theta$.
\medskip

\begin{theorem}[{\cite[Proposition 4]{gobe:makh:10a}}]
\label{thm:gm_1}
Under the condition \SDE, for $0<\theta<1$ and a bounded $g$, the following assertions are equivalent:
\begin{enumerate}[{\rm (i)}]
\item $g(X_T)\in \bes_{2,2}^\theta$.
\item $\int_0^T (T-t)^{-\theta} \E \bet \nabla_x u(t,X_t) \rag^2 dt < \infty$.
\item $\int_0^T (T-t)^{1-\theta} \E \bet D^2 u(t,X_t) \rag^2 dt < \infty$.
\end{enumerate}
\end{theorem}
\bigskip

\begin{theorem}[{\cite[Lemma 6]{gobe:makh:10a}}]
\label{thm:gm_2}
Under the condition \SDE, for $0<\theta \le1$ and a bounded $g$, the following assertions are equivalent:
\begin{enumerate}[{\rm (i)}]
\item $g(X_T)\in \bes_{2,\infty}^\theta$.
\item $\sup_{t\in [0,T)} (T-t)^{1-\theta} \E \bet \nabla_x u(t,X_t) \rag^2 < \infty$.
\item For $0<\theta<1$ we have that
      $\sup_{t\in [0,T)} (T-t)^{2-\theta} \E \bet D^2 u(t,X_t) \rag^2 < \infty$.
\end{enumerate}
\end{theorem}
\bigskip
Theorems \ref{thm:gm_1} and \ref{thm:gm_2} generalize results obtained in
\cite{geis:geis:04} and \cite{geis:hujo:07}. We see that the fractional smoothness index $\theta$ measures exactly the growth rate of the derivatives of the associated PDE solved by $u$ (see question (Q1) in the introduction). 

The two above theorems are also valid if $u$ is computed using the risk-neutral measure $\Q$ 
(i.e. $u_\Q(t,x)= \E_\Q (g(X_T) | X_t=x)$), while the other $L_2$-norms are computed under $\P$. 
For instance, for $0<\theta < 1$ the equivalence of (i) and (ii) of Theorem \ref{thm:gm_1} 
becomes $g(X_T)\in \bes_{2,2}^{\theta,\P}$ if and only if $\int_0^T (T-t)^{-\theta} \E_\P \bet \nabla_x u_\Q(t,X_t) \rag^2 dt < \infty$ 
where  we have indicated explicitly if the $L_2$-norms or conditional expectations are computed under $\P$ or $\Q$.
This property can be established following \cite{gobe:makh:10a} and the proof of \cite[Lemma 7]{gobe:makh:10b}. 
This accommodates well the fact that the price functions are usually computed under the risk-neutral measure, while hedging is 
made under the historical probability (see question (Q2) in the introduction).
\medskip

\begin{proof}[Simplified proof of Theorem \ref{thm:gm_2}]

We sketch the proof in the simple case where $X=W$ is a linear Brownian motion, $d=1$ and  
$\theta\in(0,1)$. First, $(u(t,W_t)=\E(g(W_T)|{\cal F}_t))_{t \le T}$ is a martingale in 
$L_2(\P)$. In addition, 
for any fixed $0<\delta< T$ the processes $(\nabla_x u(t,W_t))_{t \le T-\delta}$ and 
$(D^2 u(t,W_t))_{t \le  T-\delta}$ are $L_2(\P)$-martingales. This property is obtained by
checking that $\nabla_x u$ and $D^2u $ both solve the parabolic heat equation and that
certain integrability assumptions are satisfied.
Then by It\^o's formula, one obtains for $0\le s\le t<T$ that
\begin{eqnarray}
      g(W_T)-u(t,W_t) 
& = & \int_t^T \nabla_x u(s,W_s)dW_s,\label{eq:1}\\
      \nabla_x u(t,W_t)-\nabla_x u(s,W_s)
& = & \int_s^t D^2 u(r,W_r)dW_r.\label{eq:2}
\end{eqnarray}
From the It\^o isometry, one deduces from (\ref{eq:1}) that $\E|g(W_T)-u(t,W_t)|^2=\int_t^T\E| \nabla_x u(s,W_s)|^2 ds$ and it follows that  $(ii)\Rightarrow (i)$.
Similarly from (\ref{eq:2}) one obtains 
\[     \E|\nabla_x u(t,W_t)|^2
   \le 2 \E|\nabla_x u(0,W_0)|^2+2 \int_0^t \E|D^2 u(r,W_r)|^2dr \]
which proves $(iii)\Rightarrow (ii)$.
Finally, we show $(i)\Rightarrow (iii)$. Standard computations give that
\begin{eqnarray*} 
      (D^2u)(t,W_t) 
& = & D^2_z  \int_\R g(x) \frac{e^{-\frac{(x-z)^2}{2(T-t)}}}{\sqrt{2\pi(T-t)}} dx
      \bigg|_{z=W_t} \\
& = & \int_\R g(x) \frac{ (x-z)^2-(T-t) }{(T-t)^2 }
      \frac{e^{-\frac{(x-z)^2}{2(T-t)}}}{\sqrt{2\pi(T-t)}} dx\bigg|_{z=W_t} \\
& = & \E \left(g(W_T) \frac{(W_T-W_t)^2-(T-t)}{(T-t)^2}|{\cal F}_t\right) \\
& = & \E \left([g(W_T)-\E(g(W_T) |{\cal F}_t)]
      \frac{(W_T-W_t)^2-(T-t)}{(T-t)^2}|{\cal F}_t\right)
\end{eqnarray*}
which implies that
\[ \| D^2u(t,W_t)\|_{L_2(\P)} \le \frac{\|W_1^2-1\|_{L_2(\P)}}{T-t} \|g(W_T)-\E(g(W_T) |{\cal F}_t)\|_{L_2(\P)} \]
so that we are done.
\end{proof}


\section{Connection to real interpolation}
\label{sec:interpolation}

Let us connect Definition \ref{def:Besov_X} to the classical notion of
fractional smoothness which also explains the notation we have used.
In particular, this connection will make clear the difference between 
$\bes_{2,\infty}^\theta$ and $\bes_{2,2}^\theta$.

\begin{definition}[\cite{BL,BS}]
Assume a couple of Banach spaces $(E_0,E_1)$ so that $E_1$ is continuously 
embedded into $E_0$. Given $x\in E_0$ and $0<\lambda<\infty$, the $K$-functional 
is given by
\[    K(x,\lambda;E_0,E_1) 
   := \inf \{ \| x_0 \|_{E_0} + \lambda \| x_1 \|_{E_1} :
              x = x_0 + x_1 \}. \]
Moreover, given $0<\theta<1$ and $1\le q \le \infty$ we define the real interpolation norm 
\[    \| x \|_{\theta,q} 
   :=  \left \| \lambda^{-\theta}  K(x,\lambda;E_0,E_1) 
      \right \|_{L_q\left ( (0,\infty),\frac{d\lambda }{\lambda} \right )} \]
and the space
$    (E_0,E_1)_{\theta,q} 
   := \left \{ x\in E_0:  \| x \|_{\theta,q} < \infty \right \}$.
\end{definition}
With our setting ($E_1$ is continuously embedded into $E_0$) we obtain the following lexicographical 
ordering of the interpolation spaces:
\[ E_1 \subseteq (E_0,E_1)_{\theta,p}  
       \subseteq (E_0,E_1)_{\theta,q}
       \subseteq (E_0,E_1)_{  \eta,r} 
       \subseteq E_0 \]
for all $0<\eta < \theta < 1$, $1 \le p \le q \le \infty$ and all $1\le r \le \infty$.
\bigskip

We apply this concept to the analysis on the Wiener space, which needs to introduce some standard 
notation (see \cite[Sections 1.1 and 1.2]{nual:06}). Let $H$ be a separable real 
Hilbert space with the scalar product denoted by $\langle .,. \rangle_H$
and $(\mathcal{M},\Sigma,\mu)$ be a complete probability space. 
We assume an isonormal family $g=\{g_h:h\in H\}$ of centered Gaussian random variables, i.e.
\[ \E_\mu (g_h \ g_k) = \langle h,k \rangle_H \quad \mbox{for all $h,k\in H$,} \]
and that $\Sigma$ is the completed $\sigma$-field generated by the random variables $\{g_h:h\in H\}$.

For each $n\geq1$, we denote by ${\cal H}_n$ the closed linear subspace of $L_2(\mu)$ generated 
by the random variables $\{H_n(g_h):h\in H, \|h\|_H=1\}$ where 
\begin{equation}\label{eq:hermite}
H_n(x)=\frac{ (-1)^n }{\sqrt{n!}}e^{\frac{ x^2 }{2 }} \frac{ d^n }{dx^n }(e^{-\frac{ x^2 }{2 }}),
\end{equation}
i.e.  the $n$-th Hermite polynomial. ${\cal H}_0$ is the set of constants. ${\cal H}_n$ is the so-called Wiener chaos of 
order $n$ and we define by $J_n:L_2(\mu)\to L_2(\mu)$ the orthogonal projection onto ${\cal H}_n$. The following 
orthogonal decomposition is known as the Wiener chaos decomposition:  
\[ L_2(\mu)=\bigoplus_{n=0}^\infty {\cal H}_n. \]
Now, we are in a position to define the Malliavin Sobolev space and Malliavin Besov space.
\begin{definition}
The Malliavin Sobolev space $\D_{1,2}(\mu) \subseteq L_2(\mu)$ is given by
\[ \D_{1,2}(\mu) := \left \{ Z \in L_2(\mu) : 
               \| Z \|_{\D_{1,2}(\mu)} := \left ( \sum_{n=0}^\infty (n+1) \| J_n Z \|_{L_2(\mu)}^2 \right )^\frac{1}{2}
                                 < \infty \right \}.\]
Moreover, given $0<\theta<1$ and $1 \le q \le \infty$, we define the Malliavin Besov space
\[ \B_{2,q}^\theta(\mu) := (L_2(\mu),\D_{1,2}(\mu))_{\theta,q} \]
of fractional smoothness $\theta$ with fine parameter $q$.
\end{definition}
\bigskip
We use this construction in the case that $H=\ell_2^d$ and
$\mathcal{M}=\R^d$, $\Sigma$ is the completion of the Borel $\sigma$-algebra on $\R^d$ and
$\mu=\gamma_d$ is the $d$-dimensional standard Gaussian measure. The family of Gaussian
random variables is given by
\[  g_h(x) := \langle x,h \rangle
    \sptext{1}{for}{1}
    x\in \mathcal{M}=\R^d   
    \sptext{1}{and}{1}
    h \in H=\ell_2^d. \]
To make the connection between the definitions of $\bes_{2,q}^\theta$ and 
$\B_{2,q}^\theta(\gamma_d)$ for $q\in\{2,\infty\}$ we let, as before, $(W_t)_{t\in [0,1]}$ be the 
standard $d$-dimensional Brownian motion on $(\Om,¸\ftn,\P,(\ftn_t)_{t\in [0,1]})$. 
Then we have
\bigskip

\begin{theorem}[{\cite[Corollary 2.3]{geis:hujo:07}}]
\label{cor:interpolation-eta-q-convergence-r^d}
For $0<\theta < 1$, $1 \le q \le \infty$, and $g\in L_2(\gamma_d)$ one has
\[       \| g \|_{\B_{2,q}^\theta(\gamma_d)} 
  \sim_c \| g \|_{L_2(\gamma_d)} + 
         \left \| (1-t)^{-\frac{\theta}{2}} 
         \noo M_1 - M_t \rrm_{L_2(\P)}
         \right \|_{L_q([0,1), \frac{dt}{1-t})} \]
where $M_t := \E \left  (g(W_1)|\ftn_t\right )$
and $c\ge 1$ depends on $(\theta,q)$ only.
\end{theorem}
\bigskip
Applying this theorem to $q=\infty$ gives that
\[       \| g \|_{\B_{2,\infty}^\theta(\gamma_d)} 
  \sim_c \| g \|_{L_2(\gamma_d)} + 
         \sup_{0\le t \le 1} (1-t)^{-\frac{\theta}{2}} \noo M_1 - M_t \rrm_{L_2(\P)},
         \]
whereas $q=2$ gives that
\[       \| g \|_{\B_{2,2}^\theta(\gamma_d)} 
  \sim_c \| g \|_{L_2(\gamma_d)} + 
         \left ( \int_0^1 (1-t)^{-1-\theta} 
         \noo M_1 - M_t \rrm_{L_2(\P)}^2 dt \right )^\frac{1}{2}
         \]
which brings us back to Definition \ref{def:Besov_X}.
\bigskip

\paragraph{Multi-dimensional Black-Scholes-Samuelson model.}
This is a log-normal model which dynamics on the price and the log-price
can be written as
\begin{eqnarray*}
      dS^i_t
& = & S^i_t\left (\sum_{j=1}^d\sigma_{ij} dW^j_t+\mu_i dt\right ),\quad 1\le i\le d,\\
      X_t^i 
& = & \log(s^i_0)+\sum_{j=1}^d \sigma_{ij} W_t^j + (\mu_i-\frac 12\sigma_i^2)t,
\end{eqnarray*}
where $\sigma_i :=  \sqrt{\sum_j \sigma_{ij}^2}$. Assume that $(\sigma_{ij})_{i,j=1}^d$ is invertible. To the payoff function $S\mapsto h(S)$, we associate $g(x_1,...,x_d) 
   := h\left ( \left ( s_0^i e^{\sum_{j=1}^d \sigma_{ij} x_j + \mu_i - \frac{\sigma_i^2}{2} }\right )_{i=1}^d \right ) $. From this we see that 
\[ g \in \B_{2,q}^\theta(\gamma_d)
   \sptext{1.5}{if and only if}{1.5}
   h(S_1) \in \bes_{2,q}^\theta \]
for $q\in\{ 2,\infty \}$ and $g\in L_2(\gamma_d)$.

\begin{remark}\rm
In the case $\theta=1$ we get that
\[ g \in \D_{1,2}(\gamma_d) 
   \sptext{1.5}{if and only if}{1.5}
   h(S_1) \in \bes_{2,\infty}^1 \]
for all $g\in L_2(\gamma_d)$. This can be checked by using arguments from the
proof of \cite[Corollary 2.3]{geis:hujo:07}.
\end{remark}


\section{Examples}
\label{sec:examples}
In this section, we provide examples of random variables $Z=g(X_T)$ for which we determine the fractional smoothness.

\begin{example}[Lipschitz function]\rm
The case, where the fractional smoothness is obvious, is the Lipschitz case. 
Assume a Lipschitz function $g:\R^d\to\R$ with constant $L\ge 0$, i.e.
$|g(x)-g(y)| \le L |x-y|$ and assume \SDE. Then one has that
\equa
      \E \bet g(X_T) - \E       ( g(X_T)| \ftn_t ) \rag^2
&\le& \E \bet g(X_T) -            g(X_t          ) \rag^2 \\
&\le& L^2 \E \bet X_T-X_t \rag^2 \\
&\le& L^2 c^2 (T-t),
\tion
using standard estimates on the increments of $X$. Hence, $g(X_T)\in \bes_{2,\infty}^1$.
This example includes call and put payoffs, i.e.
$g(x)=(x-K)^+$ or $g(x)=(K-x)^+$.
\end{example}
\medskip
Exactly the same argument as above yields for $\theta$-H\"older functions $g$ with $\theta\in(0,1)$
that $g(X_T)\in \bes_{2,\infty}^\theta$. But the situation is here not as clear as one expects
as shown by

\begin{example}\rm
Assume the setting \GBMs and that 
\[ h_\theta(x) := (x-K)_+^\theta \]
for some $K>0$ and $0<\theta<1/2$. Then it is shown in \cite[Lemma 2]{gobe:tema:01} (under more general assumptions) 
that $\E|D^2 u(t,X_t)|^2\le c(T-t)^{-3/2+\theta}$ so that Theorem \ref{thm:gm_2} gives that
\[ h_\theta(S_T)\in\bes_{2,\infty}^{\theta+\frac 12}. \]
For $1/2<\theta<1$ one gets $h_\theta(S_T)\in\bes_{2,\infty}^1$.
\end{example}

\begin{example}[Binary option]\rm
Generally, indicator functions yield to a fractional smoothness of order $\frac 12$.
In the case $X=W$, $d=1$ and $g(x)=\mathbf{1}_{[L,\infty)}(x)$ with $L\in\R$ one 
has
\equa
u(t,x)          & = & \P(x+W_T-W_t\geq L)={\cal N}\kla \frac{x-L}{\sqrt{T-t}}\mer,\\
\nabla_x u(t,x) & = & \frac{1}{\sqrt{2\pi(T-t)}}\exp \kla -\frac{(x-L)^2}{2(T-t)}\mer,
\tion
so that
\[ \E|\nabla_x u(t,W_t)|^2 \sim_c \frac{1}{\sqrt{T-t}} \]
and $g(W_T)\in\bes_{2,\infty}^{\frac 12}$ because of Theorem \ref{thm:gm_2}.
This can be extended to the  \SDEs case as follows: Our assumption guarantees that $X$ has 
a transition density $\Gamma$ such that
\[     \Gamma (s,x;t,y) 
   \le \sqrt{\frac{\kappa}{2\pi(t-s)}} e^{- \frac{1}{2} \frac{(x-y)^2}{\kappa(t-s)}}
    =  \kappa \gamma_{\kappa(t-s)}(x-y) \]
for some $\kappa>0$ and all $0\le s < t \le T$, where $\gamma_t$ is the Gaussian density
with zero expectation and variance $t$ (see \cite{Friedman-0}). Then we can compute that
\equa
&   & \E \bet \mathbf{1}_{[L,\infty)}(X_T) - \E \kla \mathbf{1}_{[L,\infty)}(X_T)| \ftn_t \mer \rag^2 \\
&\le& \E \bet \mathbf{1}_{[L,\infty)}(X_T) - \mathbf{1}_{[L,\infty)}(X_t)\rag \\
& = & \P(X_T<L\le X_t) + \P(X_t < L \le X_T) \\
&\le& \kappa^2 [   \P(W_{\kappa T} < L-x_0 \le W_{\kappa t}) 
                 + \P(W_{\kappa t} < L-x_0 \le W_{\kappa T})] \\
&\le& c \sqrt{T-t}
\tion
where $X_0=x_0$
so that $\mathbf{1}_{[L,\infty)}(X_T) \in \bes_{2,\infty}^{\frac{1}{2}}$. 
The application in financial mathematics is done via 
$S_t = e^{X_t}$ which gives, for a positive strike $K>0$, 
\[   \mathbf{1}_{\{ S_T \ge K      \} }
   = \mathbf{1}_{\{ X_T \ge \log K \} } \in \bes_{2,\infty}^{\frac{1}{2}}. \]
In our context the fractional smoothness of jump functions (under different assumptions)
was considered in \cite{gobe:tema:01,geis:02,geis:geis:04}.
In certain multi-dimensional settings one can deduce for
$g(x)=\mathbf{1}_{\{ x_1\ge K_1, \cdots, x_d\ge K_d\}}$ (or variants of it) the same fractional smoothness from the
$1$-dimensional case.
Finally, the indicator function $g(x)=\mathbf{1}_D(x)$ of a ${\cal C}^2$-domain $D$ also leads to $g(X_T)\in\bes_{2,\infty}^{\frac 12}$ (see  \cite[Proposition 1.2]{gobe:muno:05}).
\end{example}

\begin{example}[an extreme case]\rm
By the choice of the previous examples, we emphasize that random variables $g(X_T)=h(S_T)$, usually 
used in financial applications, belong to a space $\bes_{2,\infty}^{\theta}$ for some $\theta\in(0,1]$. 
However, it is not true that $\cup_{\theta\in(0,1]}\bes_{2,\infty}^{\theta}= L_2(\P)$. The following 
result gives a way to construct $g(W_1)$ belonging to $L_2(\P)$ (here $W$ is the linear Brownian motion) but 
$g(W_1)\notin \bes_{2,\infty}^{\theta}$ for all $\theta\in(0,1]$:
\medskip

\begin{proposition}[\cite{geis:hujo:07}]
Let $0<\theta<1$, $g=\sum_{k=0}^\infty \alpha_k H_k\in L_2(\gamma_1)$, where 
$(H_k)_{k\geq 0}$ is the orthogonal basis of 
Hermite polynomials defined in  {\rm (\ref{eq:hermite})}. Then 
\[\mbox{$g(W_1)\in \bes_{2,\infty}^\theta$ \quad if and only if }\quad 
    \sup_{0\le t <1} (1-t)^{1-\theta} \sum_{k=1}^\infty k t^{k-1} \alpha_k^2 < \infty. \]
\end{proposition}
\medskip
Approximation properties as described in Section \ref{subsec:tracking_error}
for $g$ with $g(W_1)\in L_2(\P) \setminus \bigcup_{0<\theta \le 1} \bes_{2,\infty}^{\theta}$
were studied in \cite{Hujo:06} and \cite{sepp:08}.
\end{example}


\section{Applications}
\label{section:appli}
In this section we discuss some applications in stochastic finance which lead us to the
fractional smoothness as introduced above. As mentioned at the beginning, a central role is played by the tracking error that arises when 
discrete time hedging is used, instead of a continuous time strategy. For the sake of convenience, we briefly recall the notation:
\begin{itemize}
\item the option payoff at maturity $T$ is $Z=h(S_T)$;
\item the fair price function is $H(t,x)= \E_\Q (h(S_T)| S_t=x)$;
\item the $n$ re-balancing dates are defined by a deterministic time-net $\tau=(t_i)_{i=0}^n$ with $0=t_0 < t_1 < \cdots < t_{n-1} < t_n = T$;

\item the resulting tracking error process $C(Z,\tau)=(C_t(Z,\tau))_{t\in [0,T]}$ is given by
\[   C_t(Z,\tau) 
   := \E_\Q(Z|\ftn_t) - \E_\Q Z - \sum_{i=0}^{n -1}\nabla_x H(t_i,S_{t_i})\cdot (S_{t_{i+1}\wedge t}-S_{t_i\wedge t}). \]
\end{itemize}

\subsection{Weak limits of error processes}
Weak limits of stochastic processes have been intensively studied in the literature; see, 
for instance, \cite{kurt:prot:91:2,jaco:97,jaco:sh:03}. For the particular problem of the weak convergence of 
the tracking error the reader is referred to 
\cite{root:80,gobe:tema:01,haya:mykl:05,geis:toiv:09}. To formulate our results, we let 
$\widetilde{W}=(\widetilde{W}_t)_{t\ge 0}$ 
be a standard Brownian motion starting at zero defined on some auxiliary probability space, 
where we may and do assume that all paths are continuous. 
In the following $\Longrightarrow_{C[0,s]}$ stands for the weak convergence in 
$C[0,s]$ for some $s>0$.

In this paragraph we assume that $T=1$ 
\footnote{With $T=1$ we are in accordance with the quoted literature 
          that used Hermite polynomials. Of course, we could do a re-scaling to $T>0$ afterwards.}
and that $S$ is the standard geometric Brownian motion, 
i.e. the setting of {\GBM} and $\P=\Q$. The following result is the starting point of this section:
\medskip

\begin{theorem}[\cite{gobe:tema:01}]
Let $\tau_n=(i/n)_{i=0}^n$ be the equidistant time-nets and
let $Z:= \mathbf{1}_{[K,\infty)}(S_1)$ be the payoff of a digital option with 
strike price $K>0$. Then one has that
\[ \sqrt{n} C_1(Z,\tau_n) 
   \Longrightarrow
   \widetilde{W}_{\frac{1}{2}\int_0^1 
         \bet S_t^2 \frac{\partial^2 H}{\partial x^2}(t,S_t)\rag^2 dt} \]
where $\Longrightarrow$ denotes the weak convergence as $n$ goes to infinity.
\end{theorem}
\bigskip
The remarkable fact is that the weak limit is not square-integrable. In the following we describe
a way to increase the integrability of the weak limit. This is of particular interest for risk management purposes,
as a higher integrability gives better tail-estimates. The idea
is to use adapted time-nets that are more concentrated close to maturity. They are defined 
as follows: Given a parameter $\theta\in (0,1]$, we define the nets $\tau^{n,\theta}$ by
\[ t_k^{n,\theta} := 1 - \kla 1 - \frac{k}{n} \mer^\frac{1}{\theta}. \]
For $\theta=1$ we have the equidistant time-nets, i.e. $t_k^{n,1} = \frac{k}{n}$.
Now we have
\bigskip

\begin{theorem}[\cite{geis:toiv:09}]
\label{theorem:main_theorem_Besov_new}
Let $0<\theta\le 1$, $Z= h(S_1) \in L_2(\P)$ and $0 \le s < 1$. Then
\[ 
   (\sqrt{n} C_t (Z,\tau ^{n, \theta}))_{t\in [0,s]} 
   \Longrightarrow_{C[0,s]}
      \kla \widetilde{W}_{\int_0^t \frac{(1-r)^{1-\theta}}{2\theta} \bet 
           S_r^2 \frac{\partial^2 H}{\partial x^2}(r,S_r)\rag^2 dr} 
      \mer_{t\in [0,s]}. \]
Moreover, the following assertions are equivalent:
\begin{enumerate}[{\rm (i)}]
\item One has $h(S_1)  \in \bes_{2,2}^\theta$ for $0<\theta<1$ or
      $h(S_1)  \in \bes^1_{2,\infty}$ for $\theta=1$.
\item On some stochastic basis there exists a continuous 
      square-integrable martingale $M=(M_t)_{t\in [0,1]}$ such that
      $\sqrt{n} C(Z,\tau^{n,\theta}) \Longrightarrow_{C[0,1]} M$.
\item For 
      \[  A:= \int_0^1 \frac{(1-t)^{1-\theta}}{2\theta} 
              \bet S_t^2 \frac{\partial^2 H}{\partial x^2}(t,S_t)\rag^2 dt \]
      one has that $\E A< \infty$ and 
      \[                          \sqrt{n} C(Z,\tau ^{n, \theta})
         \Longrightarrow_{C[0,1]} 
         \kla \widetilde{W}_{ \mathbf{1}_{\{A<\infty\} }\int_0^t \frac{(1-r)^{1-\theta}}{2\theta}
              \bet S_r^2 \frac{\partial^2 H}{\partial x^2}(r,S_r)\rag^2 dr} 
              \mer_{t\in [0,1]}. \]
\end{enumerate} 
\end{theorem}
\bigskip

The theorem above gives us one way to consider the $L_p$-setting for 
$2\le p < \infty$. Given a differentiable function $\psi:(0,\infty)\to \R$ we let
\[ (A\psi )(x) :=   x \psi'(x) - \psi(x). \]
In the following $AH(t,x)$ means that $A$ acts on the $x$-variable of the function $H(t,x)$.

\begin{definition}\label{definiton:D_beta}\rm 
For $h(S_1)\in L_2(\P)$, $0<\theta < 1$, and $0\le t <1$ we let
\[    D^{S,\theta}_t h(S_1) 
   := \frac{1-\theta}{2}\int_0^1 (1-u)^{-\frac{1+\theta}{2}}
      [AH(u\wedge t,S_{u\wedge t}) - AH(0,S_0)] du. \]
For $\theta=1$ and  $t\in [0,1)$ we let $D^{S,1}_t h(S_1) :=  AH(t,S_t) - AH(0,S_0)$.
\end{definition}
\bigskip
 
The process $D^{S,\theta} h(S_1)=(D^\theta_t h(S_1))_{t\in [0,1)}$ is 
a quadratic integrable martingale on the half open time interval $[0,1)$.
Using the Riemann-Liouville operator of partial integration the process 
$D^{S,\theta} h(S_1)$ can be interpreted as a fractional differentiation of 
order $\theta$ in $x$ (see \cite{geis:toiv:09}).
The point of the construction of 
$D^{S,\theta} h(S_1)$ is that we may have $L_p$-singularities of
$S_t\frac{\partial H}{\partial x}(t,S_t)$ 
as $t\uparrow 1$ whereas $D^{S,\theta} h(S_1)$ remains $L_p$-bounded.
\bigskip

\begin{theorem}[\cite{geis:toiv:09}]
For $2\le p <\infty$, $0<\theta \le 1$, and $Z=h(S_1)\in L_2(\P)$
the following assertions are equivalent:
\begin{enumerate}[{\rm (i)}]
\item On some stochastic basis there exists a continuous 
      $L_p(\P)$-integrable martingale $M$ such that
      $\sqrt{n} C(Z,\tau^{n,\theta}) \Longrightarrow_{C[0,1]} M$.
\item The martingale $D^{S,\theta} h(S_1)$ is bounded in 
      $L_p(\P)$.
\end{enumerate}
\end{theorem}


\subsection{$L_2$-estimates of the tracking error}
\label{subsec:tracking_error}

In this section we work in the 1-dimensional martingale case assuming \SDEs with $\sigma(t,x)=\sigma(x)$ 
and $b(t,x)=-\frac 12 \sigma^2(x)$ (meaning $\P=\Q$). The payoff function $h$ is polynomially bounded and the option maturity 
is $T>0$. We remind the reader about the time-nets $\tau^{n,\theta}$ given by
\[ t_k^{n,\theta} := T \kla 1 - \kla 1 - \frac{k}{n} \mer^\frac{1}{\theta}\mer \]
and that for $\theta=1$ we obtain the equidistant nets.
Let us first check what quadratic hedging error one can expect at all if the portfolio 
is re-balanced $n$-times. The answer is the rate $1/\sqrt{n}$ as shown by

\begin{theorem}[{\cite[Theorem 2.5]{geis:geis:04}}]
Assume that there are no constants $c_0,c_1\in\R$ such that
$h(S_T) = c_0+c_1 S_T$ a.s. Then
\[ \inf_{n=1,2,...\atop
         0=t_0<\cdots < t_n=T}  n^{\frac 12} \| C(h(S_T),(t_k)_{k=0}^n)\|_{L_2(\P)} > 0 \]
where the infimum is taken over deterministic time-nets.
\end{theorem}

This was extended to the case of random time-nets in \cite{geis:geis:06}
in the case of the geometric Brownian motion.
\smallskip

Now we continue with the case of equidistant time-nets which are often
used in discretizations.

\paragraph{Equidistant time-nets.}
Here a starting point is the following result of Zhang:

\begin{theorem}[{\cite[Theorem 2.4.1]{zhan:99}}]\label{thm:zhang:99}
Assume that $h:\R\to \R$ is a Lipschitz function. Then we have that
\[   \lim_n n^{\frac 12} \|C(h(S_T),\tau^{n,1})\|_{L_2(\P)} \in [0,\infty). \]
\end{theorem}
\bigskip
This is the result one would expect: Given a Lipschitz payoff, the $L_2$-rate of the error
is $1/2$ for equidistant nets. But this is not the case in general as shown in
\bigskip

\begin{theorem}[{\cite[Theorem 1]{gobe:tema:01}}]
\label{thm:gobet-temam-binary}
For $h(x)=\mathbf{1}_{[K,\infty)}(x)$ for some $K>0$ we have that
\[ \lim_n n^{\frac 14} \|C(h(S_T),\tau^{n,1})\|_{L_2(\P)} \in (0,\infty). \]
\end{theorem}
\bigskip
This means that the $L_2$-approximation rate for the binary option is $n^{1/4}$ if one uses
equidistant nets. The two above results also hold true for appropriate 
$\Q\not = \P$ (i.e. $S$ is not martingale) where the outer $L_2$-norm is computed w.r.t.
the historical probability $\P$ (cf. the remarks after Theorem \ref{thm:gm_2}).
\smallskip

Theorems \ref{thm:zhang:99} and \ref{thm:gobet-temam-binary} lead
naturally to two questions: What is the reason for the 
rate $1/4$ and, secondly, can one improve the rate $1/4$? Both questions can be answered by 
the usage of the concept of fractional smoothness.

\begin{theorem}[{\cite[Theorems 2.3 and 2.8]{geis:geis:04}}]
\label{thm:Geiss-Geiss-04-1}
For $0<\theta \le 1$ and a polynomially bounded $h:(0,\infty)\to \R$ the following assertions 
are equivalent:
\begin{enumerate}[{\rm (i)}]
\item $h(S_T)\in\bes_{2,\infty}^\theta$.
\item $\sup_n n^\frac{\theta}{2} \| C(h(S_T),\tau^{n,1})\|_{L_2(\P)}<\infty$.
\end{enumerate}
\end{theorem}
\bigskip
In particular, it turns out that $h(S_T)\in \D_{1,2}$ 
if and only if
\[ \sup_n n^\frac{1}{2} \| C(h(S_T),\tau^{n,1})\|_{L_2(\P)}<\infty, \]
see \cite[Theorem 2.6]{geis:geis:04}, where $\D_{1,2}$ is the Malliavin Sobolev space
obtained from the construction in Section \ref{sec:interpolation} with
$H=L_2[0,T]$ and $g_h:= \int_0^T h(t)dW_t$.
\smallskip

For the binary option one has in Theorem \ref{thm:Geiss-Geiss-04-1} that
$\theta=1/2$ (cf. Example 3 in Section \ref{sec:examples}). This recovers the rate $1/4$ obtained in 
Theorem \ref{thm:gobet-temam-binary}.

\paragraph{Non equidistant time-nets.} Next we show how to obtain the optimal rate $n^{1/2}$ by a suitable choice of the trading dates (see question (Q4) in 
Section \ref{sec:introduction}). We can combine \cite[Lemmas 3.2 and 5.3]{geis:geis:04} and \cite[Lemma 3.8]{geis:hujo:07}
to get

\begin{theorem}\label{thm:Geiss-Geiss-Hujo}
For $0<\theta \le 1$ and a polynomially bounded $h:(0,\infty)\to \R$
the following assertions are equivalent:
\begin{enumerate}[{\rm (i)}]
\item $\int_0^T (T-t)^{1-\theta} \E \bet S_t^2 \frac{\partial^2 H}{\partial x^2}(t,S_t)\rag^2 dt < \infty$.
\item $\sup_n n^\frac{1}{2} \| C(h(S_T),\tau^{n,\theta})\|_{L_2(\P)}<\infty$.
\end{enumerate}
\end{theorem}
\bigskip
For $0<\theta < 1$ (and at least a bounded $h$) the condition 
Theorem \ref{thm:Geiss-Geiss-Hujo}(i) is equivalent to 
\begin{enumerate}
\item [(i')] $h(S_T) \in \bes_{2,2}^\theta$
\end{enumerate}
which can be checked by using Theorem \ref{thm:gm_1}.
For the binary option this gives that
\[ \sup_n n^\frac{1}{2} \| C(\mathbf{1}_{[K,\infty)}(S_T),\tau^{n,\eta})\|_{L_2(\P)}<\infty \]
for any strike $K>0$ and $0<\eta<1/2$.
\bigskip

For the next two theorems we assume that $T=1$, that $S_t=e^{W_t-\frac t2}$ and that 
$h$ might be general, i.e. not polynomially bounded. The formulation of 
Theorem \ref{thm:Geiss-Geiss-Hujo} in the language of the interpolation spaces
introduced in Section \ref{sec:interpolation} gives

\begin{theorem}[{\cite[Theorem 3.2]{geis:hujo:07}}] 
\label{theorem:interpolation-two}
For $0<\theta\le 1$ and $h(S_1) \in L_2(\P)$ the following assertions 
are equivalent:
\begin{enumerate}[{\rm (i)}]
\item $h(e^{\cdot -(1/2)}) \in \B_{2,2}^\theta (\gamma_1)$ if
      $0<\theta <1$ and $h(e^{\cdot -(1/2)}) \in \D_{1,2}(\gamma_1)$ if $\theta=1$.
\item $\sup_n n^\frac{1}{2} \| C(h(S_1),\tau^{n,\theta})\|_{L_2(\P)}<\infty$.
\end{enumerate}
\end{theorem}

And Theorem \ref{thm:Geiss-Geiss-04-1} can be extended in this context to the full
scale of real interpolation spaces as

\begin{theorem}[{\cite[Theorem 3.5]{geis:hujo:07}}]
\label{theorem:eta-q-interpolation-equidistant}
For $1\le q \le \infty$, $0<\theta<1$ and $h(S_1) \in L_2(\P)$ the following assertions 
are equivalent:
\begin{enumerate}[{\rm (i)}]
\item $h(e^{\cdot -(1/2)}) \in \B_{2,q}^\theta (\gamma_1)$.
\item $\left \| \left ( n^{\frac{\theta}{2}-\frac{1}{q}} a_n 
               \right )_{n=1}^\infty \right \|_{\ell_q} <\infty$
      for $a_n:= \| C(h(S_1),\tau^{n,1})\|_{L_2(\P)}$.
\end{enumerate}
\end{theorem}

\paragraph{Concluding remarks}

\begin{enumerate}[(i)]
\item The higher dimensional case for $X$ was considered 
      in the literature as well. Roughly speaking, one can analogously obtain upper bounds, however
      precise lower
      bounds as in the one-dimensional case are still missing. This is due to the fact that 
      a characterization of the $L_2$-error proved in \cite[Theorem 4.4]{geis:02} and
      \cite[Lemma 3.2]{geis:geis:04} is missing for higher dimensions. However, 
      after Zhang \cite{zhan:99} started with the regular case,
      Temam \cite{tema:03} extended results from \cite{gobe:tema:01} to higher dimensions 
      and Hujo \cite{Hujo:05} used non-uniform time-nets to improve the approximation rates 
      for certain irregular payoffs to the optimal rate $1/\sqrt{n}$
      in this setting.

\item Sepp\"al\"a \cite{sepp:08} found a criterion to characterize under certain conditions
      that there is a constant $c>0$ such that
      \[ \inf_{\tau = (t_i)_{i=0}^n \atop
               0=t_0<\cdots<t_n=1} \| C(h(S_1),\tau)\|_{L_2(\P)} \le \frac{c}{\sqrt{n}} \]
      where deterministic time-nets are taken.
      It should be noted that one has a {\em non-linear} approximation problem as the time-nets
      may change for fixed $n$ from payoff to payoff $h$.

\item In the above discussion, the time-nets $\tau$ are deterministic. Alternatively, one can
      allow the time-nets to be stochastic and adapted. This issue has been handled by
      \cite{mart:patr:99} using optimal stopping tools. The estimation of convergence rates is 
      an open question. However, it was shown in \cite{geis:geis:06} that the random time-nets
      do not improve the best possible approximation rate $1/\sqrt{n}$ in the case \GBMs when in the 
      $n$-th approximation a sequence of $n$ stopping times is used.

\item Similar studies can be performed when studying the Delta-Gamma hedging strategies. Instead of 
      hedging the payoff using only the asset, we use other traded options written on the same asset. 
      For a one-dimensional asset, if the price of the additional option is $(P(t,S_t))_{0\le t\le T}$, 
      the numbers of options $P$ and assets to hold at time $t_i$ are respectively equal to 
      \[ \delta^P_{t_i}:=\frac{\partial^2_S H(t_i,S_{t_i})}{\partial^2_S P(t_i,S_{t_i})}
         \sptext{.7}{and}{.7}
         \delta^S_{t_i}:= \partial_S Ht_i,S_{t_i})-\frac{\partial^2_S H(t_i,S_{t_i})}{\partial^2_S 
         P(t_i,S_{t_i})}\partial_S P(t_i,S_{t_i}). \]
      In \cite[Theorem 6]{gobe:makh:10b}, considering a multi-dimensional Black-Scholes model, it 
      is established that for an exponentially bounded payoff such that 
      $g(X_T) \in\bes_{2,\infty}^\theta$ for some $0<\theta < 1$,
      the use of equidistant time-nets leads to the same convergence rate $1/n^{\theta/2}$ as for the delta hedging 
      strategy. On the contrary, the use of non equidistant time-nets $\tau^{n,\eta}$ with $0<\eta<\theta/2$ enables us
      to obtain the improved convergence rate $1/n$.
\end{enumerate}


\section{Further developments}
\label{sct:devpt}
\subsection{Backward stochastic differential equations} 

Makhlouf and the second author applied in \cite{gobe:makh:10a} the concept of fractional 
smoothness to backward stochastic differential equations of the type
\[ Y_t = g(X_T) + \int_t^T f(s,X_s,Y_s,Z_s) ds - \int_t^T Z_s dW_s \]
where $X=(X_t)_{t\in [0,T]}$ is our forward diffusion and
the generator $f$ is continuous in its four arguments, continuously differentiable in
$(x,y,z)$ with uniformly bounded derivatives. These equations are particularly useful in stochastic finance, 
since they allow to take into account market frictions and constraints (we refer to \cite{elka:peng:quen:97} 
for a more complete account on this subject).

Solving numerically this type of equation is a challenging issue since it concerns a non-linear problem (due to the
generator $f$), generally defined in a multi-dimensional setting. One possible approach consists in approximating 
the BSDE using a discrete-time dynamic programming equation (see \cite{zhan:04,bouc:touz:04,lemo:gobe:wari:06} among
others). One of the main error contribution is related to the $L_2$-regularity on $Z$, defined by 
\[   {\cal E}(Z,\tau)
    = \sum_{i=1}^n \int_{t_{i-1}}^{t_i} \| Z_t - Z_{t_{i-1}} \|_{L_2(\P)}^2 dt. \]
If $f$ were equal to 0, then the $Z$-component is given by $z_t=\nabla_x u(t,X_t)\sigma(t,X_t)$ where
$u(t,x)=\E(g(X_T)|X_t=x)$. Studying the $L_2$-regularity of $z$ is thus very similar to the analysis 
of the tracking error presented in Section \ref{section:appli}. Additionally, using BSDE techniques, one can prove explicit upper
bounds for the difference $Z-z$
\medskip

\begin{theorem}[{\cite[Corollary 14]{gobe:makh:10a}}]
Assume \SDEs and $g(X_T)\in \bes_{2,\infty}^\theta$ for $0<\theta\le 1$. Then, for some $c>0$, one has that
\equa
       \left|Z_t-z_t\right| 
&\leq& c \int_t^T\frac{\sqrt{\E\left[\left|g(X_T)-\E( g(X_T)|{\cal F}_s)\right|^2|{\cal F}_t\right]}}{T-s}ds+c(T-t),\\
       \E\left|Z_t-z_t\right|^2 
&\leq& c (T-t)^\theta.
\tion
\end{theorem}
\medskip
Taking advantage of this approximation result close to the time singularity, we can prove that the estimate of 
${\cal E}(z,\tau)$ (linear case) transfers to ${\cal E}(Z,\tau)$ (non-linear case) and get
\medskip

\begin{theorem}[{\cite[Theorem 21]{gobe:makh:10a}}]
Assume \SDE, $g(X_T)\in \bes_{2,\infty}^\theta$ and that $0<\eta < \theta < 1$ or
$\eta=\theta=1$. Then one has that
\[{\cal E}(Z,\tau^{n,\eta}) \le \frac{c}{n}. \]
\end{theorem}
\bigskip
In \cite{geis:geis:gobe:09}, extensions of the above in different directions are discussed.

\subsection{L\'evy processes}

An extension of the results of \cite{geis:hujo:07} to L\'evy Processes
is done by C. Geiss and Laukkarinen in \cite{geis:lauk:10}. Moreover, 
Tankov and Brod\'en proved in \cite{tank:brod:09} results along the line of \cite{gobe:tema:01}.

\subsection{Multigrid Monte-Carlo Methods}

In the context of Multigrid Monte-Carlo Methods it turned out that the concept of fractional
smoothness is useful as well. The reader is referred to the papers of Avikainen
\cite{Avi_09a,Avi_09b}.


\end{document}